\documentclass[a4paper,11pt]{amsart}
\usepackage[cp1251]{inputenc}
\usepackage[english]{babel}
\usepackage{amsmath,amssymb,euscript,amsthm,amsfonts}
\usepackage{url}
\newtheorem{thm}{Theorem}
\theoremstyle{definition}
\newtheorem{defn}{Definition}
\theoremstyle{remark}

\begin{document}

\renewcommand{\proofname}{Proof}
\makeatletter \headsep 10 mm \footskip 10 mm
\renewcommand{\@evenhead}%
{\vbox{\hbox to\textwidth{\strut \centerline{{\it Mathieu Dutour}}} \hrule}}

\renewcommand{\@oddhead}%
{\vbox{\hbox to\textwidth{\strut \centerline{{\it Adjacency method for extreme Delaunay polytopes
%If the title is too long for a line, please, write our its shorten
% variant
}}} \hrule}}

%Papers for the Proceedings of the Third Vorono\"\i \ Conference of
%the Number Theory and Spatial Tessellations must be prepared in
%\LaTeX \ and according to the following scheme:

\begin{center}
{\Large\bf Adjacency method for extreme Delaunay polytopes}
\end{center}

\begin{center}
Mathieu Dutour (France) \footnote{Research financed by EC's IHRP Programme, within the Research Training Network ``Algebraic Combinatorics in Europe,'' grant HPRN-CT-2001-00272.}

\end{center}

\vskip 15pt

\begin{quote}
The {\em hypermetric cone} is defined as the cone of semimetrics satisfying the {\em hypermetric inequalities}. Every {\em Delaunay polytope} corresponds to a ray of this polyhedral cone. The Delaunay polytopes, which correspond to extreme rays are called {\em extreme}.

We use this polyhedral cone and the {\em closest vector problem} to present a new technique that allow to find, from a given extreme Delaunay polytope, some new ones.

Then, we show some examples of applications of this technique in low-dimensions.
\end{quote}

\vskip 20pt
A {\em distance vector} $(d_{ij})_{0\leq i< j\leq n}\in R^N$ with $N={n+1\choose 2}$ is called an {\em $(n+1)$-hypermetric} if it satisfies the following {\em hypermetric inequalities}:
\begin{equation*}
H(b)d=\sum_{0\leq i<j\leq n} b_ib_jd_{ij}\leq 0\mbox{~for~any~}b=(b_i)_{0\leq i\leq n}\in Z^{n+1}\mbox{~with~}\sum_{i=0}^{n}b_i=1\,\,.
\end{equation*}
The set of distance vectors, satisfying all hypermetric inequalities, is called the {\em hypermetric cone} and denoted by $HYP_{n+1}$.

In fact, $HYP_{n+1}$ is a polyhedral cone (see (7) p.~199). Lovasz (see (7) p.~201-205) gave another proof of it and bound $\max |b_i|\leq n!2^n{2n \choose n}^{-1}$ for any vector $b=(b_i)_{0\leq i\leq n-1}$ defining a facet of $HYP_n$.

There is a many-to-many correspondence between Delaunay polytopes with
$p$ vertices and the elements of the hypermetric cone $HYP_{p}$ (see
(7)). So, the study of Delaunay
polytopes is equivalent, in a way, to the study of hypermetrics.
The Delaunay polytopes, whose corresponding face is an extreme ray,
are called {\em extreme}.

Two extreme ray of a polyhedral cone are called {\em adjacent} if they generate
a two-dimensional face of this cone. Therefore, it induces a natural adjacency
relation between extreme Delaunay polytopes. Our main purpose here is to
explain how, given an extreme Delaunay polytope, one can find the adjacent
extreme Delaunay polytopes.

This technique allow us to find some new extreme Delaunay polytopes and is
illustrated in some low-dimensional cases.

The following {\em closest vector problem} is a classic problem of discrete geometry and is heavily used in our computations:

{\it Closest Vector Problem:

\nopagebreak
(i) Given a lattice $L\subset R^n$, a vector $x\in R^n$ and $r>0$, test if there exist $v\in L$ such that $\Vert v-x\Vert<r$.

(ii) Given a lattice $L\subset R^n$, a vector $x\in R^n$ and $r>0$, find all the vectors $v\in L$ such that $\Vert v-x\Vert=r$.
}

This classic problem is NP-hard (see (8)). There are several software packages doing this computation (see (11) and (14)).

%Moreover, the dimension of faces of $HYP_{p}$ allows to define the
%notion of {\em rank} of a Delaunay polytope (see \cite{DGL92}).

\section{The hypermetric cone}
For more details on the material of this Section see Chapters $13$--$16$ of (7). For the use of affine basis and hypermetrics for computing combinatorial types of Delaunay polytopes, see (9).\\

%\cite{DGL92} show that the hypermetric on any generating subset of an extreme Delaunay polytope (see above) lies on an extreme ray of $HYP_n$ and that a hypermetric, lying on an extreme ray of $HYP_n$, is the square of Euclidean distance on a generating subset of extreme Delaunay polytope of dimension at most $n-1$.
%On every set $A=\{ v_1, \dots, v_m\}$ of vertices of a Delaunay polytope $P$, one can define a distance function $d_{ij}=\Vert v_i-v_j\Vert^2$. The function $d$ turns out to be a metric and, moreover, a hypermetric. On the other hand, every hypermetric is the square of Euclidean distance on a generating set of vertices of a Delaunay polytope of a lattice.

Denote by $S(c,r)$ the sphere of center $c$ and radius $r$.
For every family ${\bf A}=\{ v_0, \dots, v_m\}$ of vertices of a Delaunay polytope $P$ circumscribed by the sphere $S(c,r)$ one can define a distance function $d_{\bf A}$ by $(d_{\bf A})_{ij}=\Vert v_i-v_j\Vert^2$. The function $d_{\bf A}$ turns out to be a hypermetric by the following formula (see (1) and (7) p.~195) :
\begin{equation*}
H(b)d_{\bf A}=\sum_{0\leq i,j\leq m} b_ib_j(d_{\bf A})_{ij}=2(r^2-\Vert \sum_{i=0}^m b_iv_i - c\Vert^2)\leq 0\,.
\end{equation*}
On the other hand, Assouad has shown in (1) that {\em every} $d\in HYP_{n+1}$ can be expressed as $d_{\bf A}$ with ${\bf A}$ being a family of vertices of a Delaunay polytope $P$ of dimension {\em less or equal} to $n$.

\begin{defn}
Let $P$ be an $n$-dimensional Delaunay polytope with vertex-set ${\bf V}$.

(i) A family  $v_0$, \dots, $v_n$ of vertices of $P$ is called an {\em affine basis} if for all $v\in {\bf V}$ there exists an unique vector $b=(b_i)_{0\leq i\leq n}\in Z^{n+1}$, such that
\begin{equation*}
\sum_{i=0}^{n}b_i v_i=v\mbox{~and~}\sum_{i=0}^{n}b_i=1\;.
\end{equation*}

(ii) The Delaunay polytope $P$ is called {\em basic} if it has at least one affine basis. The vertices of an affine basis are called {\em basic vertices}.

\end{defn}

All Delaunay polytopes of dimension less, than $6$, are basic (see (9)).
The adjacency method, described later, applies only to basic Delaunay
polytopes. There is no proof that all Delaunay polytopes are basic and the
only indication for this conjecture is that all known Delaunay polytopes are
basic (see (7)).

The {\em metric cone} $MET_{n+1}$ is defined as the set of functions $d=(d_{ij})_{0\leq i<j\leq n}$ satisfying all triangular inequalities $d_{ij}\leq d_{ik}+d_{jk}$; such inequalities correspond to the hypermetric inequality with vector $b$, such that $b_i=b_j=1$, $b_k=-1$ and $b_l=0$, otherwise. Therefore $HYP_{n+1}\subset MET_{n+1}$.

Given $d\in HYP_{n+1}$ we define 
\begin{equation*}
Ann(d)=\{b\in Z^{n+1}\mbox{~:~}\sum_{i=0}^n b_i=1\mbox{~and~}H(b)d=0\}.
\end{equation*}

\begin{thm}\label{theoremDescription}
Let $P$ be an $n$-dimensional extreme Delaunay polytope and ${\bf B}=\{v_0,\dots, v_n\}$ an affine basis of $P$. Then the following statements hold:

(i) The mapping $b\rightarrow \sum_{i=0}^n b_i v_i$ establishes a bijection between $Ann(d)$ and the vertex-set of $P$.

(ii) It is possible to compute (in finite time) the list of hypermetric inequalities $H(b_i)d\leq 0$ with $1\leq i\leq m$ incident to $d_{\bf B}$. 

(iii) $P$ has at least $\frac{(n+2)(n+1)}{2}-1$ vertices. If it has exactly $\frac{(n+2)(n+1)}{2}-1$ vertices, then all $b_i$ correspond to facets of $HYP_{n+1}$. If $P$ has more vertices, then it is possible to compute in finite time which of the $b_i$ correspond to facets of $HYP_{n+1}$.
\end{thm}
\begin{proof}
Let $S(c,r)$ be the empty sphere around $P$ and let us write the affine basis as $v_0$, \dots, $v_n$. 
The equality $H(b)d_{\bf B}=0$ corresponds to 
\begin{equation*}
\begin{array}{rcl}
H(b)d_{\bf B}&=&2(r^2-\Vert \sum_{i=0}^n b_iv_i - c\Vert^2)=0\\
&=&2(r^2-\Vert v_0-c+\sum_{i=1}^n b_i(v_i-v_0)\Vert^2)=0\,.
\end{array}
\end{equation*}
This establishes a bijection between solutions of the equation $H(b)d_{\bf B}=0$ and vertices of the Delaunay polytope $P$.

Moreover, the solution-set is exactly a closest vector problem;
therefore, it is solvable in finite time. Note that all $e_i=(0,\dots, 1,\dots, 0)$ with $0\leq i\leq n$ belong to $Ann(d_{\bf B})$; for those vector $H(e_i)=0$, therefore we do not list them in the list of vector $(b_i)_{1\leq i\leq m}$.

The Delaunay polytope is extreme; therefore, the rank of the system of equations $H(b_i)d=0$ is equal to $\frac{n(n+1)}{2}-1$. This implies $m\geq \frac{n(n+1)}{2}-1$. Since the vertex-set contains also the affine basis, it has at least $n+1+m=\frac{(n+1)(n+2)}{2}-1$ vertices.

Any extreme ray of $HYP_{n+1}$ is incident to at least $\frac{n(n+1)}{2}-1$ facets. Therefore, if $P$ has $\frac{(n+1)(n+2)}{2}-1$ vertices, then all the $H(b_i)$ are facets of $HYP_{n+1}$.

Let us assume that a facet $H(b_i)$ is redundant. Then, $H(b_i)$ can be expressed as a sum with positive coefficients of other facets of $HYP_{n+1}$:
\begin{equation*}
H(b_i)=\sum_{j\in J} \alpha_j H(b'_j)\mbox{~with~}\alpha_j>0
\end{equation*}
The relation $H(b_i)d_{\bf B}=0$ implies that for all $j\in J$, we have $H(b'_j)d_{\bf B}=0$. Therefore, $b'_j$ belongs to the list of $b_i$. So, the problem of redundancy is reduced to the one for the list of facets $H(b_i)$ with $1\leq i\leq m$, which is solvable by linear programming (see polyhedral FAQ\footnote{\url{http://www.ifor.math.ethz.ch/~fukuda/polyfaq/polyfaq.html}}).

\end{proof}
Above theorem is very useful for finding hypermetric facets. The Delaunay polytopes, whose number of vertices is exactly $\frac{(n+1)(n+2)}{2}-1$ are especially useful for the adjacency method, by the simplicity of the computation of their adjacencies.
The Schl\"afli polytope and the infinite serie of extreme Delaunay polytopes in (12) satisfy this lower bound.

\begin{thm}
$HYP_8$ has at least $298592$ facets and at least $86$ orbits of facets.
\end{thm}
\begin{proof}
We know that all facets $H(b)$ of $HYP_n$ yield facets $H(b')$ of $HYP_{n+1}$, with $b'$ obtained by adding a zero to $b$. Since $HYP_7$ has $14$ orbits of facets (see (10)), this yields $14$ orbits.

Moreover, all hypermetric facets of $CUT_n$ are facets of $HYP_n$. In (4), a conjectural list of orbits of facets of $CUT_8$ is computed. This yields $16$ orbits of facets of $HYP_8$. 

We know two extreme Delaunay polytopes in dimension $7$: the Gosset polytope and a polytope found by Erdahl and Rybnikov (see (12)). By computing all affine basis of them and using above theorem we obtain some more facets. It is known (see (7) p.~229) that the switchting of a facet of $HYP_n$ is again a facet of $HYP_n$. Combining all this we obtain $298592$ facets of $HYP_8$ in $86$ orbits.

\end{proof}

\begin{defn}
(i) Two Delaunay polytopes $P_1$ and $P_2$ are called {\em isomorphic} if there exists an isometry transforming one into the other.

(ii) The {\em automorphism group} $Aut(P)$ of a Delaunay polytope $P$ is the set of all isometries leaving it invariant.
\end{defn}
Denote by ${\bf V}_1$ and ${\bf V}_2$ the vertex-sets of $P_1$ and $P_2$, respectively.
In order to test the existence of such an isometry, it suffices to test if there exist a mapping $\phi:{\bf V}_1\rightarrow {\bf V}_2$ satisfying to $d(\phi(x),\phi(y))=d(x,y)$ for all $x,y\in {\bf V}_1$. The test of existence of an isometry is therefore a combinatorial problem. If the set of possible pairwise distances is $\{d_1,\dots, d_h\}$, then the isomorphy problem becomes a problem of isomorphy of association schemes. The problem of computing the automorphism group of a Delaunay polytope is the same as of computing the automorphism group of an association scheme.
Unfortunately, we do not know about a program doing such computations.
%Propose it

Therefore, it will be presented here some algorithms doing computations using the nauty program (see (13)). Problem (ii) can be solved easily by associating, to every possible distance $d$ between two vertices of $P$ the graph $G_d$, formed by making two vertices adjacent if $d(x,y)=d$. So, we have
\begin{equation*}
Aut(P)=\cap_{i=1}^{h} G_{d_i},
\end{equation*}
Problem (i) can be solved by associating to $P_1$ and $P_2$ a graph encoding pairwise distances, but this approach leads to graphs, whose size is too big for the nauty program. So, we use procedures, which are not guaranteed to work in any cases but which have worked so far in all cases considered of extreme Delaunay polytopes.

Take two Delaunay polytopes $P_1$ and $P_2$ and compute their skeleton graphs (i.e. the graph formed by their vertices with two vertices adjacent if they generate a two dimensional face), $G_1$ and $G_2$, respectively.
If those graphs are not isomorphic, then $P_1$ and $P_2$ are not isomorphic. If they are isomorphic, then there exists an isomorphism $\phi:G_1\rightarrow G_2$. If this mapping satisfies $d(\phi(x),\phi(y))=d(x,y)$, then $P_1$ and $P_2$ are isomorphic. This approach does not always work, i.e. some automorphisms of $G_1$ do not correspond to isometries of $P_1$ and therefore the mapping $\phi$ computed is not the right one.

A similar approach is to consider the {\em ridge graph} (i.e. the graph formed by their facets with two facets being adjacent if their intersection is of dimension $n-2$) $G'_1$ and $G'_2$ of $P_1$ and $P_2$. If $G'_1$ and $G'_2$ are not isomorphic, then $P_1$ and $P_2$ are not isomorphic. If they are isomorphic, then there exists an isomorphism $\phi':G'_1\rightarrow G'_2$. Since every vertex is an intersection of some facets, this mapping lifts to an isomorphism $\phi:G_1\rightarrow G_2$. We then test if $\phi$ defines an isometry. In all cases considered, $\phi$ was an isometry, i.e. we were able to check the isomorphy of extreme Delaunay polytopes efficiently. There is no reason to think that this will always be the case.

\section{The Adjacency Method}

We consider in this section the details of the adjacency algorithm that takes an extreme ray $e$ of $HYP_{n+1}$, corresponding to an $n$-dimensional extreme Delaunay polytope and finds the adjacent extreme rays $(e_i)_{1\leq i\leq m}$ in $HYP_{n+1}$.

%The idea is the following: given an $n$-dimensional extreme Delaunay polytope, we want to find a corresponding extreme ray $e\in HYP_{n+1}$, we then find the adjacent extreme ray and this gives us some new extreme rays and so, some new extreme Delaunay polytopes.

The problem in doing this computation is that the description of facets of $HYP_n$ is known only for $n\leq 7$ (see (9), (10) and (3)). 
Anyway, using the complete list of facets is not a good idea, since, for example, $HYP_8$ has at least $298,592$ facets.

The algorithm, which we use, is as follows:
\begin{enumerate}
\item Given an initial extreme ray $e$, find the list of hypermetric vectors $(b_i)_{1\leq i\leq m}$ such that $H(b_i)e=0$. We then take as initial list ${\bf F}$ of facets all vectors $(b_i)_{1\leq i\leq m}$ plus all permutations of $(1^2, -1, 0^{n-2})$ (they correspond to triangular inequalities).
\item Define the cone ${\bf C}({\bf F})$ by taking all hypermetric inequality in ${\bf F}$. 
\item Since $e$ is an extreme ray of the cone ${\bf C}({\bf F})$, we can find the extreme rays $(e_j)_{1\leq j\leq p}$ of ${\bf C}({\bf F})$ that are adjacent to $e$.
\item Test if $e_j$ is hypermetric, using the closest vector problem.
\item If the ray $e_j$ is hypermetric, then we finish. If not, then some hypermetric inequalities $H(b)d\leq 0$ are violated. So, we add those $b$ to ${\bf F}$ and go back to step 2.
\end{enumerate}
Since the cone $HYP_{n+1}$ is polyhedral, the inner loop will eventually finish.

The key step in above algorithm is the ability to check if $H(b)d\leq 0$ is true for all vector $b\in Z^{n+1}$ with $\sum_i b_i=1$. This is equivalent to
\begin{equation*}
H(b)d_{\bf B}=2(r^2-\Vert v_0-c+\sum_{i=1}^n b_i(v_i-v_0)\Vert^2)\leq 0\,.
\end{equation*}
If we denote by $L$ the lattice generated by the family $(v_i-v_0)_{1\leq i\leq n}$, then the problem is expressed as
\begin{center}
Does there exist $v\in L$ such that $\Vert v_0-c+v\Vert< r$?
\end{center}
This problem is a closest vector problem (i) except that the scalar product (which comes from the distance $d$) is not always positive definite.

We associate to $d$ the Gram matrix $G=(g_{ij})$ defined by $g_{ij}=\frac{1}{2}(d_{i0}+d_{j0}-d_{ij})$.

If $G$ admits a negative eigenvalue, then one can find a vector $v\in L$, such that $\Vert v_0-c+v\Vert <0$, thereby solving the problem.

If $G$ is positive but not positive definite, then this means that $L$ is of dimension inferior to $n$ and that the $n+1$ vectors $v_i$ form a Delaunay polytope of a lower dimensional lattice, which is not a simplex. In that case the algorithm is as follows:
\begin{enumerate}
\item Find the rank $r$ of the matrix $G$.
\item Find a family ${\bf F}=\{v_{i_0}, \dots, v_{i_r}\}$, such that all points $v_i$ can be expressed in terms of ${\bf F}$ with integer coefficient (this is not always possible).
\item If preceding step has succeeded, then the problem takes the form
\begin{center}
Does there exist $b_{i_j}$, such that $\Vert \sum_{j=0}^r b_{i_j}v_{i_j}-c\Vert\leq r$?
\end{center}
The corresponding Gram matrix is positive definite. Therefore, this is a closest vector problem of type (i).

\end{enumerate}

%Associated to $d$ is a $n\times n$ Gram matrix $G=(g_{ij})$ defined by $g_{ij}=\frac{1}{2}(d_{i0}+d_{j0}-d_{ij})$. It is well known (see (7) p.????) that $G$ is positive. Let us assume first that $G$ is positive definite. Then $G$ defines an Euclidean distance on the lattice $Z^n$; moreover, the $n+1$ points $v_0=(0,\dots, 0)$ and $v_i=(0,\dots,  1, \dots, 0)$ with $1\leq i\leq n$ are vertices of a Delaunay polytope $P$ of $Z^n$ for this scalar product. The problem is now to find all vertices of $P$ and our method is as follow:
%\begin{enumerate}
%\item Find the center $c$ and the radius $r$ of the empty ellipsoid $S(c,r)$ containing the $n+1$ points $v_i$ with $0\leq i\leq n$.
%\item Finding all points of $P$ is then the same as solving the equation
%\begin{equation*}
%{}^{t}(x-c) G(x-c)=r^2
%\end{equation*}
%But this problem is closest vector problem (ii).
%\end{enumerate}

\section{Example of application of the Adjacency Method}

In (6) (see also (7)) an extreme $15$-dimensional Delaunay polytope
with $135$ vertices, i.e. having the minimal number of vertices, is given.
We applied the adjacency method to one of its affine basis and found one new
extreme Delaunay polytope, let us denote it by $ED_8$ of dimension $8$.
This polytope belongs to an infinite serie $(ED_n)_{n\geq 6}$ of extreme
Delaunay polytopes (see (10)); this serie was found from $ED_8$.

We then apply the adjacency decomposition method to the polytope $ED_8$ and found $24$ extreme Delaunay polytope of dimension $8$.

It is well-known (see (7)) that the Gosset polytope $Gos$ is an extreme $7$-dimensional Delaunay polytope. All orbits of affine basis of $Gos$ were found by direct enumeration.

We apply the Adjacency Method to all orbits of affine basis of $Gos$ and found an extreme Delaunay polytope of dimension $7$ with $35$ vertices and a symmetry group of size $1440$. In fact, this extreme Delaunay polytope was already found by Erdahl and Rybnikov in (12).
We conjecture that there is no other extreme Delaunay polytopes in dimension $7$.

\vskip 10 pt

\centerline{\Large REFERENCES} \vskip 7 pt

\begin{enumerate}
\small
\item P. Assouad (1982), {\em Sous-espaces de $L^1$ et in\'egalit\'es hyperm\'etriques}, Compte Rendus de l'Acad\'emie des Sciences de Paris, {\bf 294(A)} 439--442.

\item E.P. Baranovski, {\em Simplexes of $L$-subdivisions of Euclidean spaces}, Mathematical Notes, {\bf 10} (1971) 827-834.

\item E.P. Baranovskii, {\em The conditions for a simplex of $6$-dimensional lattice to be $L$-simplex}, (in Russian) Nauchnyie Trudi Ivanovo state university. Mathematica, {\bf 2} (1999) 18--24.

\item T. Christof and G. Reinelt (2001), {\em Decomposition and parallelization techniques for enumerating the facets of combinatorial polytopes}, Internat. J. Comput. Geom. Appl., {\bf 11-4} 423--437.

\item M. Deza and M. Dutour (2003), {\em The hypermetric cone on seven vertices}, to appear in Experimental Mathematics.

\item M. Deza, V.P. Grishukhin, and M. Laurent, {\em Extreme hypermetrics and L-polytopes}, in G.Hal\'asz et al. eds {\em Sets, Graphs and Numbers, Budapest (Hungary), 1991}, {\bf 60} {\em Colloquia Mathematica Societatis J\'anos Bolyai}, (1992) 157--209.

\item M. Deza and M. Laurent (1997), {\em Geometry of Cuts and Metrics}, Berlin, Heidelberg, New York: Springer Verlag.

\item I. Dinur, G. Kindler, S. Safra (1998), {\em Approximating CVP to within almost-polynomial factors is NP-hard}, 39th Annual IEEE Symposium on Foundations of Computer Science.

\item M. Dutour (2003), {\em The six-dimensional Delaunay polytopes}, to appear in European Journal of Combinatorics.

\item M. Dutour (2003), {\em Infinite serie of extreme Delaunay polytopes}, to appear in European Journal of Combinatorics.

\item M. Dutour (2003), {\em Lattice-CVP}, \url{http://www.liga.ens.fr/~dutour/CVP/index.html}

\item R. Erdahl and K. Rybnikov, {\em Supertopes}, \url{http://faculty.uml.edu/krybnikov/}

\item B. McKay, {\em The nauty program}, \url{http://cs.anu.edu.au/people/bdm/nauty/}

\item F. Vallentin (1999), {\em ShVec}, \url{http://www-m10.ma.tum.de/~vallenti/shvec.02nov99.tar.gz}

\end{enumerate}

\vskip 20 pt

\noindent {\it Hebrew University Jerusalem and \'Ecole Normale Sup\'erieure Paris}

\vskip 2 pt

\noindent {\it e-mail: {\tt Mathieu.Dutour@ens.fr}}

\end{document}